\documentclass[a4paper,11pt]{article}
\usepackage{amsmath}
\usepackage{amssymb}
\usepackage{amsfonts}
\usepackage{graphicx}
\usepackage{color}
\usepackage{comment}
\usepackage{bm}
\usepackage[round]{natbib}

\newtheorem{theorem}{Theorem}
\newtheorem{lemma}[theorem]{Lemma}

\newtheorem{remark}{\it Remark}

\newtheorem{proposition}[theorem]{Proposition}

\newcommand{\qed}{\hfill $\square$}

\title{Identifiability of asymmetric circular and cylindrical distributions}
\author{Yoichi Miyata \thanks{Faculty of Economics, Takasaki City University of Economics, 1300 Kaminamie, Takasaki, Gunma, 370-0801, Japan. E-mail: ymiyatagbt@tcue.ac.jp, ORCID:0000-0002-9822-690X}\and 
        Takayuki Shiohama \thanks{Department of Data Science, Nanzan University, 18 Yamazato-cho Showa-ku, Nagoya, 446-8673, Japan.}\and Toshihiro Abe \thanks{Faculty of Economics, Hosei University, 4342 Aihara, Machida, Tokyo, 194-0298, Japan.}
}

\begin{document}
\maketitle

\begin{abstract}
Identifiability of statistical models is a fundamental and essential condition that is required to prove the consistency of maximum likelihood estimators. The identifiability of the skew families of distributions on the circle and cylinder for estimating model parameters has not been fully investigated in the literature. In this paper, a new method combining the trigonometric moments and the simultaneous Diophantine approximation is proposed to prove the identifiability of asymmetric circular and cylindrical distributions. Using this method, we prove the identifiability of general sine-skewed circular distributions, including the sine-skewed von Mises and sine-skewed wrapped Cauchy distributions, and that of a M\"{o}bius transformed cardioid distribution, which can be regarded as asymmetric distributions on the unit circle. In addition, we prove the identifiability of two cylindrical distributions wherein both marginal distributions of a circular random variable are the sine-skewed wrapped Cauchy distribution, and conditional distributions of a random variable on the non-negative real line given the circular random variable are a Weibull distribution and a generalized Pareto-type distribution, respectively.
\end{abstract}

\section{Introduction}\label{Intro}
\hspace{1em}Identifiability is a prerequisite and essential condition for most common statistical inferences, especially for deriving the consistency of estimators. A family of statistical models is said to be identifiable if no two sets of parameter values yield the same distribution. Identifiability of statistical models has been considered in various contexts; for details we refer the reader to \cite{Ro71} and the monograph of \cite{Ra92}.

Traditional statistical models and methods are based on data taking values in Euclidean space, which relies on the linear structure of the space and the fundamental operations of linear algebra. Recently, various statistical data analyses have been introduced for manifolds because the linear statistical approaches are inappropriate for these non-linear spaces. Examples of statistical manifolds include the circle, sphere, and cylinder. For further reading on statistical data analysis in various manifolds, 
see the recent monographs of \cite{Ch12}, \cite{BB08}, \cite{DK15}, and \cite{LV17b}.

Asymmetric circular and spherical distributions have been proposed by \cite{UJ09}, \cite{AP11}, \cite{JP12}, and various others. Cylindrical models with skewed marginal distributions on the circle were proposed by \cite{AL17} and \cite{imoto2019cylindrical}. 
Because the distributions proposed by them have a skewness parameter, 
we can provide a more meaningful interpretation than that obtained by finite mixture modeling of symmetric circular distributions when the data show a clear asymmetric pattern. Alternatively, there are approaches that utilize transformations of scale \citep{JP12} and the M\"{o}bius transformation \citep{KJ10}. Although the identifiability of these skew families of the distributions on circle and cylinder is of primary importance for estimating model parameters, it has not been fully investigated in the literature. 

Because almost all finite mixture models are not identifiable due to label switching and potentially overfitting \citep[e.g., see][]{Fr06}, the generic identifiability of \cite{YS68} is discussed instead. 
Generic identifiability implies that if two finite mixture distributions in a parametric family are equal, then identical parameter vectors can be chosen, corresponding to these two distributions under suitable permutation and aggregation of component densities. For details, see \cite{Fr06}, \cite{YS68} and \citet[Definition 3.1.1]{TSM85}. To demonstrate the generic identifiability of finite mixture models, \citet[Theorem 2]{Te63} uses one-to-one transforms $m_{j}(t)$ $(j=1,2)$ of the distribution functions of any two distinct components indexed by $t$, such as the moment generating function and the characteristic function, and verifies that the ratio $m_{1}(t)/m_{2}(t)$ converges to zero or infinity as $t\to \infty$. 
\cite{HMS04} employed Teicher's approach to demonstrate the generic identifiability of finite mixtures of symmetric circular distributions including the von Mises and wrapped Cauchy distributions.
The point here is that any component density satisfying the conditions of \cite{Te63} in a finite mixture model is identifiable in the usual sense. Accordingly, the result of \cite{HMS04} indicates that the von Mises and wrapped Cauchy distributions are identifiable. 
In contrast, Teicher's approach cannot be directly applied to asymmetric circular models such as the sine-skewed distributions because, when the transforms $m_{j}(t)$ $(j=1,2)$ are set to be characteristic functions of two different models, the ratio $m_{1}(t)/m_{2}(t)$ does not generally converge when $t\to\infty$. 

In this paper, applying Teicher's approach for the trigonometric moments of circular random variables and their functions, we show the ordinary identifiability of some asymmetric circular and cylindrical distributions. 
In general, a sequence of ratios of the $p$th trigonometric moments $(p=1,2,...)$ does not have a limit as $p\to\infty$, because it includes sine and cosine functions such as $\sin (p\mu)$ and $\cos(p\mu )$ where $\mu$ is any constant in the interval $[0,2\pi )$. However, there exists a subsequence $\{ p_{n}\}_{n\in\mathbb{N}}$ such that the trigonometric functions described above converge and $p_{n}\to\infty$ as $n\to \infty$. 
For example, suppose $\mu=0.7$, and consider a subsequence $\{ p_{n} \}$ with $p_{1}=47$, $p_{2}=82$, $p_{3}=180$, $p_{4}=386,\cdots $. Then a sequence $\{p_{n}\mu \,(\textrm{mod}\, 2\pi ) \}$ is given by $1.48$, $0.85$, $0.34$, $0.02, \cdots$, where these numbers are rounded to two decimal places, and the symbol ``$x \, (\textrm{mod}\, 2\pi )$" stands for the remainder after a number $x$ is divided by $2\pi$. Hence, $\sin(p_{1}\mu )>\sin(p_{2}\mu )>\sin(p_{3}\mu )>\sin(p_{4}\mu )\geq 0$ holds. In Lemma \ref{lemma1}, we use the simultaneous Diophantine approximation to show that there exists a subsequence $\{ p_{n}\}_{n\in\mathbb{N}}$ such that a sequence of the ratios of the $p_{n}$th trigonometric moments converges and $p_{n}\to\infty$ as $n \to\infty$. The simultaneous Diophantine approximation is used to evaluate the accuracy when a vector of real numbers was approximated by a vector of rational numbers. See \cite{Sc91,Sc96} for details.

The remainder of the paper is organized as follows. In Section~\ref{sec:Main}, we propose a new approach to demonstrate identifiability whose conditions are weaker than those of \cite{Te63}. In Section~\ref{sec:applications}, we elucidate the identifiability of general sine-skewed circular distributions, including the sine-skewed von Mises (SSvM) and the sine-skewed wrapped Cauchy (SSWC) distributions, and a M\"{o}bius transformed distribution. 
In Section \ref{sec:applications02}, we demonstrate the identifiability of the cylindrical distribution of \cite{AL17}, which combines the sine-skewed von Mises distribution on the circle and the Weibull distribution on the non-negative linear line, and the cylindrical distribution of \cite{imoto2019cylindrical}, wherein the generalized Pareto distribution is considered for the linear part. Section~\ref{sec:concluding} concludes the paper, and all proofs are given in Appendix~\ref{proofs}. We give fundamental trigonometric moment properties for sine-skewed circular distributions in Appendix~\ref{appendix:A}.

\section{Main results}\label{sec:Main}
\hspace*{1em}This section presents a general approach to demonstrate the identifiability for a family of the cumulative distribution functions (CDFs).
Let $F(\bm{x} |\bm{\gamma})$ be the CDF of an observed random vector $\bm{X}\in\mathbb{R}^{m}$, and let $\Gamma$ be a parameter space in $k$-dimensional Euclidean space $\mathbb{R}^{k}$.
The random vector $\bm{X}$ does not have to be a circular random variable. First, we rigorously define the identifiability of a family $\mathcal{F}:=\{ F(\bm{x} |\bm{\gamma})|\bm{\gamma}\in\Gamma\}$. 
We say that the family $\mathcal{F}$ is identifiable if for any elements $\bm{\gamma}_{1}$, $\bm{\gamma}_{2}$ in $\Gamma$, $\bm{\gamma}_{1}\ne\bm{\gamma}_{2}$ implies $F(\bm{x}|\bm{\gamma}_{1})\ne F(\bm{x}|\bm{\gamma}_{2})$. 
If the family $\mathcal{F}$ has a density function $f(\bm{x}|\bm{\gamma})$, then a family $\{ f(\bm{x} |\bm{\gamma})|\bm{\gamma}\in\Gamma\}$ of densities and $\mathcal{F}$ can be considered identical. 
Before describing our proposed approach, we must pay special attention to the non-equal condition for parameter vectors $\bm{\gamma}_{1}$ and $\bm{\gamma}_{2}$. Then, we can divide $\bm{\gamma}_{1}\ne\bm{\gamma}_{2}$ into some cases.
For example, we consider the case for a model parameter dimension of $k=2$. Then, the non-equal symbol is decomposed as follows: 
for any $\bm{\gamma}_{1}=(\gamma_{1,1},\gamma_{1,2})^{T}$ and $\bm{\gamma}_{2}=(\gamma_{2,1},\gamma_{2,2})^{T}$, 
\begin{equation}
\bm{\gamma}_{1}\ne \bm{\gamma}_{2}\qquad 
\textrm{if and only if} \qquad
\begin{cases}
\textrm{(i) } \gamma_{1,1}\ne\gamma_{2,1},\\
\qquad \textrm{or} \\
\textrm{(ii) } \gamma_{1,1}=\gamma_{2,1}\quad \textrm{and}\quad \gamma_{1,2}\ne\gamma_{2,2}.
\end{cases}
\end{equation}
We refer to (i) as ``Step 1" and (ii) as ``Step 2". 

Now, to demonstrate the identifiability of the model, we consider the following conditions, which are weaker than those of \cite{Te63}. 
Suppose that $d$ is the total number of steps, and $\phi_{i} (\bm{p}|\bm{\gamma})$ is a transform of $F(\bm{x} |\bm{\gamma})$ in Step $i$. Note that $k$ and $d$ do not necessarily have to be the same.
An example of such a transform $\phi_{i}  (\bm{p}|\bm{\gamma})$ is the moment generating function $M(\bm{p}|\bm{\gamma}):=E_{\bm{\gamma}}[\exp(\bm{p}^{T}\bm{X})]$ where $\bm{a}^{T}$ denotes the transpose of a matrix or a vector $\bm{a}$ and $E_{\bm{\gamma}}[\bullet ]$ indicates the expectation under the CDF $F(\bm{x} |\bm{\gamma})$. 
Let ${\cal M}_{\bm{p}}(F(\bm{x} |\bm{\gamma}))=(\phi_{1} (\bm{p}|\bm{\gamma}),...,\phi_{d} (\bm{p}|\bm{\gamma}))^{T}$ be a mapping with an index $\bm{p}\in\mathbb{R}^{m}$ from $\mathcal{F}$ to $\mathcal{G}:=\{ (\phi_{1} (\bm{p}|\bm{\gamma}),...,\phi_{d} (\bm{p}|\bm{\gamma}))^{T}| \bm{\gamma}\in\Gamma\}$. 
$S_{i}(\bm{\gamma}):=\{ \bm{p}\in\mathbb{R}^{m}|\phi_{i}(\bm{p}|\bm{\gamma})\textrm{ is well-defined and }$\\ $|\phi_{i}(\bm{p}|\bm{\gamma})|<\infty\}$ denotes the domain of the function $\phi_{i}(\bm{p}|\bm{\gamma})$ from $\mathbb{R}^{m}$ to $\mathbb{C}$ for each $\bm{\gamma}$ and Step $i$. $\bar{S}$ denotes the closure of set $S$ under a suitable metric $d(x,y)$. We often set $S=\mathbb{Z}$ with metric $d(x,y)=|\tan^{-1}x-\tan^{-1}y|$. Then, $\bar{S}=\mathbb{Z}\cup \{ \pm\infty\}$ becomes the extended rational number set. 
We now present the main theorem of this paper.
\begin{theorem}\label{main_theorem}
Suppose that for any $\bm{\gamma}_{1}$ and $\bm{\gamma}_{2}$ in $\Gamma$ with $\bm{\gamma}_{1}\ne \bm{\gamma}_{2}$, there exist a mapping $\mathcal{M}_{\bm{p}}:\mathcal{F}\to \mathcal{G}$, a Step $i\in \{ 1,...,d\}$, and a non-empty set $D_{i}$ such that 
\begin{list}{}{} 
\item[(a)] $D_{i}\subseteq S_{i}(\bm{\gamma}_{1})\cap S_{i}(\bm{\gamma}_{2})$, and $D_{i}$ does not depend on $\bm{\gamma}_{1}$ and $\bm{\gamma}_{2}$,
\item[(b)] there exist a point $\bm{p}_{\infty}^{i}\in \overline{D_{i}}$, a constant $c_{i}\ne 1$ and a sequence $\{ \bm{p}_{n}^{i}\}_{n\in \mathbb{N}}$ in $D_{i}$ with $\lim_{n\to\infty}\bm{p}_{n}^{i}=\bm{p}_{\infty}^{i}$ such that 
\begin{equation}
\frac{\phi_{i}(\bm{p}_{n}^{i}|\bm{\gamma}_{2})}{\phi_{i} 
(\bm{p}_{n}^{i}|\bm{\gamma}_{1})}\rightarrow c_{i}, \qquad (n\to\infty ), \label{cond(b)}
\end{equation}
\end{list}
Then, the family $\{ F(\bm{x}|\bm{\gamma})|\bm{\gamma}\in\Gamma \}$ is identifiable.
\end{theorem}
Note that the mapping $\mathcal{M}_{\bm{p}}:\mathcal{F}\to \mathcal{G}$ is not necessarily one-to-one, and among the transforms $\phi_{i}(\bm{p}|\bm{\gamma})$ $(i=1,...,d)$, some may be the same. 
In (\ref{cond(b)}), the limit is taken along some sequence $\{ \bm{p}_{n}^{i}\}_{n\in\mathbb{N}}$ whose element is in the domain $D_{i}$ of $\bm{p}$. 

\begin{remark}\label{remark1}
We consider a special case of the condition of Theorem \ref{main_theorem} that $\phi_{1}(\bm{p}|\bm{\gamma})=\cdots =\phi_{d}(\bm{p}|\bm{\gamma})$, $c_{1}=\cdots =c_{d}=0$, $\bm{p}_{n}^{1}=\cdots =\bm{p}_{n}^{d}=\bm{p}$, $\bm{p}_{\infty}^{1}=\cdots =\bm{p}_{\infty}^{d}$, $D_{1}=\cdots =D_{d}$, the mapping $\mathcal{M}_{\bm{p}}:\mathcal{F}\to \mathcal{G}$ is linear and one-to-one, and equation (\ref{cond(b)}) holds as $||\bm{p}||\to\infty$, where $||\bullet||$ is the Euclidean norm. Furthermore, the symbol $\ne$ is replaced by a total ordering $\prec$. An example of the ordering is presented in \citet[p.40]{TSM85}. Then, these conditions are equivalent to those of Theorem 2 of \cite{Te63}. 
Now, we replace the condition that equation (\ref{cond(b)}) holds  as $||\bm{p}||\to\infty$ with a weaker condition that (\ref{cond(b)}) holds for some sequence $\{ \bm{p}_{n}\}$. Even then, the generic identifiability immediately holds from the same argument as in the proof of Theorem 2 of \cite{Te63} under the relaxed conditions, which leads to an extension of the theorem. However, because this is beyond the scope of this paper, we focus instead on proving the ordinary identifiability of statistical models on manifolds.
\end{remark}

\subsection{Application to circular distributions}\label{sec:applications}
\hspace*{1em}In this section, we apply Theorem \ref{main_theorem} to asymmetric circular distributions to show their identifiability. 
Suppose that $f_{0}(\theta |\bm{\psi} )$ is a base symmetric continuous circular density of a random variable $\Theta$ whose location is set at $0$ with a parameter vector $\bm{\psi}$, that is, $f_{0}(2\pi -\theta |\bm{\psi} )=f_{0}(\theta |\bm{\psi} )$ for any $\theta \in (0,2\pi )$. 
We denote the $p$th cosine and sine moments of $\Theta$ with density $f_{0}(\theta |\bm{\psi} )$ by $\alpha_{0,p}(\bm{\psi} ):=E_{0,\bm{\psi}}\left\{ \cos (p\Theta )\right\}$ and $\beta_{0,p}(\bm{\psi} ):=E_{0,\bm{\psi}}\left\{ \sin (p\Theta )\right\}$, respectively. 
Here, we consider that a random variable $\Theta$ has the sine-skewed circular density
\begin{align}
f(\theta |\bm{\gamma})&=f_{0}(\theta -\mu |\bm{\psi} )\left\{ 1+\lambda \sin (\theta -\mu )\right\}  \label{sscd}
\end{align}
where $\bm{\gamma}=(\mu, \bm{\psi}^{T}, \lambda )^{T}$ is an element of a parameter space $\Gamma :=\{ \bm{\gamma}|\mu \in [0,2\pi ), \bm{\psi} \in \Psi, \lambda \in [-1,1]\}$, and $\Psi$ is a subset in $\mathbb{R}^{k-2}$ $(k\geq 3)$. $\mu$ is a circular location parameter. The parameter $\lambda$ controls the skewness of the distribution (\ref{sscd}), and $\lambda =0$ leads to a symmetric distribution with density $f_{0}(\theta -\mu |\bm{\psi} )$. The parameter $\bm{\psi}$ plays a role in the other shape such as unimodality and/or concentration. 
To apply Theorem \ref{main_theorem}, we employ the $p$th cosine and sine moments defined in $p\in\mathbb{Z}$, $\alpha_{p}(\bm{\gamma}):=E_{\bm{\gamma}}[\cos (p\Theta )]$ and $\beta_{p}(\bm{\gamma}):=E_{\bm{\gamma}}[\sin (p\Theta )]$, and their functions such as the $p$th mean resultant length $\rho_{p}(\bm{\gamma}) := \sqrt{\alpha_{p}(\bm{\gamma})^2 +\beta_{p}(\bm{\gamma})^2  }$ as choices for the transform $\phi_{i}(\bm{p}|\bm{\gamma})$. By combining the result of \citet[p.686]{AP11} and Appendix \ref{appendix:A}, $\alpha_{p}(\bm{\gamma})$ and $\beta_{p}(\bm{\gamma})$ are given by 
\begin{align}
\alpha_{p}(\bm{\gamma})&=\cos (p\mu )\alpha_{0,p}(\bm{\psi} )-\sin (p\mu )\lambda \{ \alpha_{0,p-1}(\bm{\psi} )-\alpha_{0,p+1}(\bm{\psi} )\} /2,\notag \\
\intertext{and}
\beta_{p}(\bm{\gamma})&=\sin (p\mu )\alpha_{0,p}(\bm{\psi} )+\cos (p\mu )\lambda \{ \alpha_{0,p-1}(\bm{\psi} )-\alpha_{0,p+1}(\bm{\psi} )\} /2 .\label{sine_cosine_moments}
\end{align}
Thereby, the $p$th mean result length becomes 
\begin{align}
\rho_{p}(\bm{\gamma})=\sqrt{\alpha_{0,p}(\bm{\psi} )^{2}+\frac{\lambda^{2}}{4}\left\{ \alpha_{0,p-1}(\bm{\psi} )-\alpha_{0,p+1}(\bm{\psi} )\right\}^{2}}. \label{sscd_mlr}
\end{align}
Let $\bm{\gamma}_{1}\ne \bm{\gamma}_{2}$ be different parameter vectors. The sequences of ratios $\alpha_{p}(\bm{\gamma}_{1})/\alpha_{p}(\bm{\gamma}_{2})$ and $\beta_{p}(\bm{\gamma}_{1})/\beta_{p}(\bm{\gamma}_{2})$ generally do not have limits as $p\to\infty$. 
However, by using Lemma \ref{lemma1} with the simultaneous Diophantine approximation, we can find a sequence $\{ p_{n}\}$ such that $p_{n}\to\infty$ as $n\to\infty$ and each of sequences of the ratios $\alpha_{p_{n}}(\bm{\gamma}_{1})/\alpha_{p_{n}}(\bm{\gamma}_{2})$ and $\beta_{p_{n}}(\bm{\gamma}_{1})/\beta_{p_{n}}(\bm{\gamma}_{2})$ converges to a non-unity constant as $n\to\infty$. 
From this fact and Theorem \ref{main_theorem}, the following result holds. The proof is given in Section \ref{proofs}.
\begin{theorem}\label{main_theorem2}
Assume that
\begin{list}{}{}
\item[(i)] for any $\bm{\psi} \in \Psi$, $\alpha_{0,1}(\bm{\psi} )\ne 0$, 
\item[(ii)] for any $\bm{\psi} \in \Psi$, there exists a constant $M=M(\bm{\psi} )>0$ such that 
\begin{align}
&\inf_{p\in\mathbb{N}}\left| \frac{\alpha_{0,p-1}(\bm{\psi} )-\alpha_{0,p+1}(\bm{\psi} )}{\alpha_{0,p}(\bm{\psi} )}\right| >\frac{1}{M}, \notag
\end{align}
\item[(iii)] there exists a constant $c\geq 0$ such that for any $\bm{\psi} \in \Psi$,
\begin{align}
\frac{\alpha_{0,p-1}(\bm{\psi} )-\alpha_{0,p+1}(\bm{\psi} )}{\alpha_{0,p}(\bm{\psi} )} =O(p^{c}), \label{cond(iii)}
\end{align}
and for any $\bm{\psi}_{1}\ne\bm{\psi}_{2}$ $(\bm{\psi}_{1},\bm{\psi}_{2}\in \Psi)$, it holds either 
\begin{align}
&\frac{\alpha_{0,p}(\bm{\psi}_{1})}{\alpha_{0,p}(\bm{\psi}_{2})}p^{c}\rightarrow 0\quad (p\to \infty )\quad \textrm{or}\quad\frac{\alpha_{0,p}(\bm{\psi}_{1})}{\alpha_{0,p}(\bm{\psi}_{2})}p^{-c}\rightarrow \infty\quad (p\to \infty ). \label{cond(iii)-2}
\end{align}
\end{list}
Then, the family $\mathcal{F}_{1}:=\{ f(\theta |\bm{\gamma})|\bm{\gamma}\in \Gamma \}$ of densities (\ref{sscd}) is identifiable.
\end{theorem}
Condition (i) is standard in symmetric circular distributions. Condition (ii) is needed to verify the identifiability of the density (\ref{sscd}) under $\lambda_{1}\ne\lambda_{2}$. Condition (iii) is a sufficient condition to ensure that $\rho_{p}(\bm{\gamma}_{1})\ne \rho_{p}(\bm{\gamma}_{2})$ holds for different parameters $\bm{\gamma}_{1}$ and $\bm{\gamma}_{2}$ with $\bm{\psi}_{1}\ne\bm{\psi}_{2}$. 
Because $p$ is an integer, the limits (\ref{cond(iii)-2}) are taken along $p\in \mathbb{Z}$ and $p\to\infty$. 
Note that conditions (i)--(iii) are imposed for the cosine moments of the symmetric base density, and can therefore be checked relatively easily. 
\begin{remark}\label{remark2}
Let $\Gamma^{*} :=\{ \bm{\gamma}|\mu \in [-\pi ,\pi ), \bm{\psi} \in \Psi, \lambda \in [-1,1]\}$ be the parameter space $\Gamma$ with $\mu \in [0,2\pi )$ replaced by $\mu \in [-\pi ,\pi )$. Then, there exists a bijective mapping between the families $\mathcal{F}_{1}$ and $\mathcal{F}_{1}^{*}:=\{ f(\theta |\bm{\gamma})|\bm{\gamma}\in \Gamma^{*} \}$. Therefore, if $\mathcal{F}_{1}$ is identifiable, then $\mathcal{F}_{1}^{*}$ is identifiable, and vice versa.
\end{remark}

As an example of (\ref{sscd}), we consider the SSWC distribution with density
\begin{equation}
f_{\textsc{SSWC}}(\theta |\mu ,\rho ,\lambda )=\frac{1 - \rho^2}{2\pi \{ 1 + \rho^2 - 2 \rho \cos (\theta  - \mu )\}}\left\{ 1+\lambda \sin (\theta -\mu )\right\} ,\label{SSWC}
\end{equation}
and suppose that $\bm{\gamma}_{\textsc{WC}}=(\mu ,\rho ,\lambda )^{T}$ belongs to a parameter space 
\begin{align}
\Gamma_{\textsc{WC}}=\left\{ \bm{\gamma}_{\textsc{WC}}\bigl| 0\leq \mu <2\pi ,\: 0<\rho <1,\: -1\leq \lambda \leq 1\right\} .\notag
\end{align} 
By combining \citet[p.696]{AP11} and Appendix \ref{appendix:A}, the $p$th cosine and sine moments are given by $\alpha_{p}(\bm{\gamma}_{\textsc{WC}}):=E\{ \cos (p\Theta )\} =\cos (p\mu )\rho^{|p|}-\sin (p\mu )\lambda (\rho^{|p-1|}-\rho^{|p+1|})/2$ and $\beta_{p}(\bm{\gamma}_{\textsc{WC}}):=E\{ \sin (p\Theta )\} =\sin (p\mu )\rho ^{|p|}+\cos (p\mu) \lambda (\rho^{|p-1|}-\rho^{|p+1|})/2$. 
Hence the $p$th mean resultant length becomes $\rho_{p}(\rho ,\lambda )=\sqrt{\rho^{2|p|}+\lambda^{2}(\rho^{|p-1|}-\rho^{|p+1|})^{2}/4}$. 

Then, the following proposition holds from Theorem \ref{main_theorem2}.
\begin{proposition}\label{thm1}
The family $\{ f_{\textsc{SSWC}}(\theta |\mu ,\rho ,\lambda )|\, (\mu ,\rho ,\lambda )^{T}\in \Gamma_{\textsc{WC}}\}$ of the SSWC distributions is identifiable. 
\end{proposition}
As another example of (\ref{sscd}), we consider the SSvM distribution with density
\begin{align}
f_{\textsc{SSvM}}(\theta |\mu ,\kappa ,\lambda )&=\frac{1}{2\pi I_{0}(\kappa )}\exp\left\{ \kappa \cos (\theta -\mu )\right\} \left\{ 1+\lambda \sin (\theta -\mu )\right\} ,\label{SSvM}
\end{align}
where $\bm{\gamma}_{\textsc{vM}}=(\mu, \kappa, \lambda)^{T}$ belongs to the parameter space 
\begin{align}
\Gamma_{\textsc{vM}}=\left\{ (\mu ,\kappa ,\lambda )\bigl| 0\leq \mu <2\pi , \kappa >0, -1\leq \lambda \leq 1\right\}
\label{eq:gamma2}
\end{align}
and $I_{\nu}(\kappa ):=(1/2\pi )\int_{0}^{2\pi}\cos (\nu \theta )\exp\{ \kappa\cos \theta\} d\theta$ is the modified Bessel function of the first kind and order $\nu \in \mathbb{Z}$. 
By combining \citet[p.691]{AP11} and Appendix \ref{appendix:A}, the $p$th cosine and sine moments are given by 
\begin{align*}
\alpha_{p}(\bm{\gamma}_{\textsc{vM}})&=\left\{ \cos (p\mu )-\frac{p\lambda}{\kappa}\sin (p\mu )\right\} \frac{I_{p}(\kappa )}{I_{0}(\kappa )}
\intertext{and}
\beta_{p}(\bm{\gamma}_{\textsc{vM}})&=\left\{ \sin (p\mu )+\frac{p\lambda}{\kappa}\cos (p\mu )\right\} \frac{I_{p}(\kappa )}{I_{0}(\kappa )},
\end{align*}
which leads to the $p$th mean resultant length $\rho_{p}(\kappa ,\lambda )=I_{p}(\kappa )\sqrt{\kappa^{2}+(p\lambda )^{2}}/(\kappa I_{0}(\kappa ))$. Then, the following proposition holds from Theorem \ref{main_theorem2}.

\begin{proposition}\label{thm2}
The family $\{ f_{\textsc{SSvM}}(\theta |\mu, \kappa, \lambda )|\, (\mu, \kappa, \lambda)^{T}\in \Gamma_{\textsc{vM}}\}$ of the SSvM distributions is identifiable.
\end{proposition}

\cite{WS12} proposed another type of asymmetric circular distribution of $\Theta$, which is not the sine-skewed circular distribution. This distribution is derived by applying the M\"{o}bius transformation to the cardioid distribution, whose density function is given by
\begin{align}
f_{\textsc{WS}}(\theta |\bm{\gamma}_{\textsc{WS}})=\frac{(1-\rho_{\alpha}^{2})h(\bm{\gamma}_{\textsc{WS}})}{2\pi \{ 1+\rho_{\alpha}^{2}-2\rho_{\alpha}\cos (\theta -\mu ) \}}, \label{MC01}
\end{align}
where
\begin{align*}
h(\bm{\gamma}_{\textsc{WS}})=1+2\bar{\rho} \left\{ \frac{\cos (\theta -\xi -\mu )-2\rho_{\alpha}\cos \xi +\rho_{\alpha}^{2}\cos (\theta +\xi -\mu )}{1+\rho_{\alpha}^{2}-2\rho_{\alpha}\cos (\theta -\mu )}\right\},
\end{align*}
and $\bm{\gamma}_{\textsc{WS}}=(\mu ,\rho_{\alpha},\bar{\rho},\xi )^{T}$ is an element of a parameter space 
\begin{align}
\Gamma_{\textsc{WS}}=\{ (\mu ,\rho_{\alpha},\bar{\rho},\xi )| 0 \leq \mu <2\pi , 0 \leq \xi <2\pi , 0<\rho_{\alpha}<1, \bar{\rho}>0\} .\label{MC_params}
\end{align}

Although the distribution (\ref{MC01}) is a special case of the family of distributions having five parameters given in Section 8 of \cite{KJ10}, it can still capture asymmetry and bimodality of distributions. 

By the results of \citet[p.607]{WS12}, the $p$th cosine and sine moments $\alpha_{p}(\bm{\gamma}_{\textsc{WS}}):=E_{\bm{\gamma}_{\textsc{WS}}}\{ \cos (p\Theta ) \}$ and $\beta_{p}(\bm{\gamma}_{\textsc{WS}}):=E_{\bm{\gamma}_{\textsc{WS}}}\{ \sin (p\Theta ) \}$ are 
\begin{align}
\alpha_{p}(\bm{\gamma}_{\textsc{WS}})&=p\bar{\rho}\rho_{\alpha}^{p-1}(1-\rho_{\alpha}^{2})\cos (p\mu +\xi )+\rho_{\alpha}^{p}\cos (p\mu )\notag
\intertext{and}
\beta_{p}(\bm{\gamma}_{\textsc{WS}})&=p\bar{\rho}\rho_{\alpha}^{p-1}(1-\rho_{\alpha}^{2})\sin (p\mu +\xi )+\rho_{\alpha}^{p}\sin (p\mu ), \quad (p\in\mathbb{Z})
\end{align}
Thus, the $p$th mean resultant length is given by 
\begin{align}
\rho_{p}(\bm{\gamma}_{\textsc{WS}})=\sqrt{p^{2}\bar{\rho}^{2}\rho_{\alpha}^{2(p-1)}(1-\rho_{\alpha}^{2})^{2}+\rho_{\alpha}^{2p}+2p\bar{\rho}\rho_{\alpha}^{2(p-1)+1}(1-\rho_{\alpha}^{2})\cos \xi } .\notag
\end{align}
Then, the following result holds.
\begin{proposition}\label{thm6}
The family $\{ f_{\textsc{WS}}(\theta |\bm{\gamma}_{\textsc{WS}}) | \bm{\gamma}_{\textsc{WS}}\in \Gamma_{\textsc{WS}}\}$ is identifiable. 
\end{proposition}
Note that although distribution (\ref{MC01}) with $\bar{\rho}=0$ reduces to the wrapped Cauchy distribution, it loses the identifiability for $\xi$ when $\bar{\rho}=0$.

\subsection{Application to cylindrical distributions}\label{sec:applications02}
\hspace*{1em}In this section, we apply Theorem \ref{main_theorem} and the result in Section \ref{sec:applications} to joint distributions of the pair $(\Theta ,X)$ where $\Theta$ is a random variable on $[-\pi ,\pi )$ corresponding to the unit circle $\mathbb{S}^{1}$ and $X$ is a random variable on the non-negative real line $\mathbb{R}^{+}$. 
First, we consider the cylindrical distribution proposed by \cite{AL17} with density
\begin{align}
f_{\textsc{AL}}(\theta ,x|\bm{\eta})=&\frac{\alpha \beta^{\alpha}}{2\pi \cosh (\kappa )}\left\{ 1+\lambda \sin (\theta -\mu )\right\} x^{\alpha -1}\notag\\
 &\qquad \times\exp \left[ -(\beta x)^{\alpha}\left\{ 1-\tanh (\kappa )\cos (\theta -\mu )\right\}\right] , \label{Abe-Ley}
\end{align}
where $\bm{\eta}=(\alpha ,\beta ,\mu ,\kappa, \lambda )^{T}$ is a parameter vector with $\alpha >0$, $\beta >0$, $-\pi \leq \mu < \pi$, $\kappa >0$ and $-1\leq \lambda \leq 1$. 
This model has some desirable properties in which the marginal distributions of $\Theta$ and $X$, the conditional distribution of $\Theta$ given $X$, and the conditional distribution of $X$ given $\Theta$ can be expressed in explicit forms. For example, the marginal distribution of $\Theta$ becomes the SSWC distribution, and the conditional distribution of $X$ given $\Theta =\theta$ becomes the Weibull distribution
\begin{align*}
f_{\textsc{AL}}(x|\theta ;\alpha ,\beta ,\kappa )&=\alpha \left[ \beta \{ 1-\tanh (\kappa)\cos(\theta-\mu)\} ^{1/\alpha}\right]^{\alpha}x^{\alpha -1}\\
 &\qquad\times\exp\left[ -\left\{ \beta \{ 1-\tanh (\kappa)\cos(\theta-\mu)\} ^{1/\alpha}x\right\}^{\alpha} \right] .
\end{align*}

We omit the details but refer the reader to \cite{AL17}.
Then, the following result holds from Theorem \ref{main_theorem} and Proposition \ref{thm2}. 
\begin{proposition}\label{thm4}
The family $\{ f_{\textsc{AL}}(\theta ,x|\bm{\eta})|\, \bm{\eta}\in \Gamma_{\textsc{AL}}\}$ is identifiable, where $\Gamma_{\textsc{AL}}:=\Gamma_{\textsc{vM}}\times \{ (\alpha ,\beta )|\alpha >0,\beta >0\}$ and $\Gamma_{\textsc{vM}}$ is defined in equation (\ref{eq:gamma2}).

\end{proposition}
Finally, as an extension of the \cite{AL17} distribution, we consider a sine-skewed generalized Pareto-type cylindrical distribution \citep{imoto2019cylindrical} with density 
\begin{align}
&f (\theta, x | \bm{\eta})
=
\left\{ 1+\lambda \sin (\theta -\mu )\right\} f_{0}(\theta ,x|\bm{\eta}), \label{gpar}
\end{align}
where 
\begin{align*}
 f_{0}(\theta ,x|\bm{\eta})=\frac{\sqrt{1-\kappa^2}}{2\pi\sigma \delta} \left(\frac{x}{\sigma}\right)^{1/\delta-1}
\left[1 + \frac{\tau}{\delta} \left(\frac{x}{\sigma}\right)^{1/\delta} \{1 - \kappa \cos(\theta - \mu) \}\right]^{-(\delta/\tau + 1)},
\end{align*}
and $\bm{\eta}=(\sigma, \delta, \tau, \mu ,\kappa, \lambda )^{T}$ is a parameter vector with $\sigma >0$, $\delta >0$, $\tau >0$, $-\pi \leq \mu < \pi$, $\kappa >0$ and $-1\leq \lambda \leq 1$. When a random vector $(\Theta ,X)$ has the density \eqref{gpar}, the marginal density of $\Theta$ is written as
\begin{align}
f(\theta |\bm{\eta})=\left\{ 1+\lambda \sin (\theta -\mu )\right\} \int_{0}^{\infty}f_{0}(\theta ,x|\bm{\eta})dx. \label{gpar_marginal}
\end{align}
Because, by equation (4) of \cite{imoto2019cylindrical}, the above integral becomes the density of the wrapped Cauchy distribution, equation \eqref{gpar_marginal} becomes the density of the SSWC distribution. 
From this, the conditional density of $X$ given $\Theta =\theta$ is 
\begin{align}
&f_{X|\Theta}(x|\theta ;\bm{\eta})=\frac{f_{0}(\theta ,x|\bm{\eta})}{\int_{0}^{\infty}f_{0}(\theta ,x|\bm{\eta})dx}\notag \\
=& \frac{1}{\sigma \delta} \left(\frac{x}{\sigma}\right)^{1/\delta -1} \{ 1 - \kappa \cos(\theta - \mu) \} \left[ 1 + \frac{\tau}{\delta} \left( \frac{x}{\sigma}\right)^{1/\delta} \{ 1 - \kappa \cos (\theta - \mu ) \}\right]^{-(\delta /\tau + 1)}.\label{gpar_cond}
\end{align}
As a result, we see that the above density is the same as equation (5) of \cite{imoto2019cylindrical}.
Let us define the parameter space of the sine-skewed generalized Pareto-type cylindrical distribution as
\begin{align*}
\Gamma_{\textsc{ISA}}= \Gamma_{\textsc{vM}} \times  
\left\{ \left.
( \delta, \tau, \sigma ) \right| \delta >0, \tau >0, \sigma> 0 
\right\}.
\end{align*}
Recall that $\Gamma_{\textsc{vM}}$ is defined in equation (\ref{eq:gamma2}). Then, we have the following result.

\begin{proposition}\label{thm5}
The family $\{ f (\theta, x | \bm{\eta}) |\, \bm{\eta}\in \Gamma_{\textsc{ISA}}\}$ of the sine-skewed generalized Pareto-type cylindrical distributions defined in (\ref{gpar}) is identifiable.
\end{proposition}

\section{Concluding Remarks}\label{sec:concluding}
\hspace*{1em} In this paper, we proposed a method to prove the identifiability of circular and cylindrical distributions under weaker conditions than those of \citet[Theorem 2]{Te63}. 
However, identification of the inverse Batschelet distribution, which is an alternative asymmetric circular distribution proposed by \cite{JP12}, remains a challenging topic. Although it is not clear at this time whether this model is identifiable or not, its characteristic function was explicitly presented by \cite{Ab15}. Therefore, there is a possibility that application of Theorem \ref{main_theorem} with this characteristic function solves this problem. 

Because Theorem \ref{main_theorem} is relatively general, it could be applied to spherical models with general skewing functions other than the sine-skewing function, which are presented in \cite{LV17a}, and toroidal models such as \cite{WS12} on $\mathbb{S}^{1}\times\mathbb{S}^{1}$. Further theoretical and methodological developments concerning asymmetric distributions on such manifolds are
expected in the near future.

Recently, \cite{MSA19} discussed the ordinary identifiability of a finite mixture of SSWC distributions when the true number of components is known. However, less is known about the generic identifiability of a finite mixture of asymmetric circular distributions. As a further topic, it would be interesting to investigate whether our approach mentioned in Remark \ref{remark1} is applicable to finite mixtures of asymmetric circular distributions and those of cylindrical distributions to demonstrate generic identifiability. Such an investigation would be valuable.
\section*{Acknowledgements}
Yoichi Miyata was supported in part by JSPS KAKENHI Grant Number 19K11863 and the competitive research expenses of Takasaki City University of Economics. Takayuki Shiohama was supported in part by JSPS KAKENHI Grant Number 18K01706 and Nanzan University Pache Research Subsidy I-A-2 for the 2021 academic year. Toshihiro Abe was supported in part by JSPS KAKENHI Grant Numbers 19K11869 and 19KK0287.

\bibliographystyle{spbasic}      


%
%

\appendix
\section{Proofs}\label{proofs}
{\sc Proof of Theorem} \ref{main_theorem}. We prove this by contradiction. Assume conditions (a) and (b), and that there exist parameters $\bm{\gamma}_{1}$, $\bm{\gamma}_{2}$ such that $\bm{\gamma}_{1}\ne\bm{\gamma}_{2}$ and $F(\bm{x}|\bm{\gamma}_{1})=F(\bm{x}|\bm{\gamma}_{2})$. 
By these assumptions, it holds that $\phi_{i}(\bm{p}|\bm{\gamma}_{2})/\phi_{i}(\bm{p}|\bm{\gamma}_{1})=1$ for every $i$ and any $\bm{p}$, which contradicts condition (b). \qed

\begin{lemma}\label{lemma1}
For any $s\in\mathbb{N}$, any $a_{1},...,a_{s}\in [0,2\pi )$ and any $\epsilon >0$, there exist infinitely many $p\in\mathbb{N}$ such that $|pa_{i} \textrm{\upshape (mod 2$\pi$)} |<\epsilon$ $(i=1,...,s)$, where ``\textrm{\upshape mod}" indicates the modulo operation.
\end{lemma}

{\sc Proof}. Fix any $s\in\mathbb{N}$ and any $\epsilon >0$. For any $a_{i}\in [0,2\pi )$ $(i=1,...,s)$, there exist constants $c_{i}\in [0,2)$ such that $a_{i}=c_{i}\pi$ $(i=1,...,s)$. 
From the simultaneous Diophantine approximation (e.g., see Theorem 1A on page 27 of \cite{Sc96}), for any natural number $Q$ with $Q\geq 2$, there exist integers $q$, $p_{1},...,p_{s}$ such that $q\in \{ 1,...,Q^{s}-1\}$, and 
\begin{equation*} 
\left| c_{i}-\frac{p_{i}}{q}\right| \leq \frac{1}{qQ},\quad (i=1,...,s)
\end{equation*}
which leads to
\begin{equation} 
p\left| c_{i}\pi-\frac{p_{i}}{q}\pi\right| \leq \frac{p\pi}{qQ}, \label{Diophantine02}
\end{equation}
where $p$ is any natural number. For simplicity, we let $c_{i}^{*}=p_{i}/q$. 
Then, choosing the number $Q$ with $2\pi Q^{-1/2}<\epsilon$ and $Q^{1/2}\in\mathbb{N}$, and putting $p=2qQ^{1/2}$ in (\ref{Diophantine02}) yields
\begin{equation} 
p\left| c_{i}\pi-c_{i}^{*}\pi\right| \leq \frac{2\pi}{Q^{1/2}}<\epsilon \qquad (i=1,...,s).  \label{Diophantine03}
\end{equation}
In addition, we have $pc_{i}^{*}\pi =2Q^{1/2}p_{i}\pi =0$ (mod 2$\pi$). 
Therefore, because $pc_{i}\pi =(Q^{1/2}p_{i})(2\pi )+(pc_{i}\pi -pc_{i}^{*}\pi )$ and $Q^{1/2}p_{i}\in\mathbb{Z}$, it follows from the inequality (\ref{Diophantine03}) that
\begin{align*}
|pc_{i}\pi \textrm{(mod 2$\pi$)}| &=|pc_{i}\pi -pc_{i}^{*}\pi | <\epsilon ,\notag 
\end{align*}
which completes the proof. Note that we can take infinitely many $p=2qQ^{1/2}$ by choosing $Q$ successfully because there are infinitely many $Q$ with $2\pi Q^{-1/2}<\epsilon$ and $Q^{1/2}\in\mathbb{N}$, and $q \geq 1$. \qed\vspace{1ex}

\hspace*{-1em}{\sc Proof of Theorem} \ref{main_theorem2}. For different parameters $\bm{\gamma}_{1}$ and $\bm{\gamma}_{2}$, we consider the following three steps:
Step 1. $\bm{\psi}_{1}\ne\bm{\psi}_{2}$, Step 2. $\bm{\psi}_{1}=\bm{\psi}_{2}$, $\lambda_{1}\ne\lambda_{2}$, Step 3. $\bm{\psi}_{1}=\bm{\psi}_{2}$, $\lambda_{1}=\lambda_{2}$, and $\mu_{1}\ne\mu_{2}$. Then, under each step, we verify the conditions in Theorem \ref{main_theorem}.

Step 1. We consider two parameter vectors $\bm{\gamma}_{1}$ and $\bm{\gamma}_{2}$ with $\bm{\psi}_{1}\ne\bm{\psi}_{2}$, and set $\phi_{1}(p|\bm{\gamma})=\rho_{p}(\bm{\gamma})^{2}$ where $\rho_{p}(\bm{\gamma})$ is defined in (\ref{sscd_mlr}). 
Then, the ratio of $\phi_{1}(p|\bm{\gamma}_{1})$ to $\phi_{1}(p|\bm{\gamma}_{2})$ is 
\begin{align}
\frac{\phi_{1}(p|\bm{\gamma}_{1})}{\phi_{1}(p|\bm{\gamma}_{2})}&=\frac{\alpha_{0,p}(\bm{\psi}_{1} )^{2}+(\lambda_{1}^{2}/4)\left\{ \alpha_{0,p-1}(\bm{\psi}_{1} )-\alpha_{0,p+1}(\bm{\psi}_{1} )\right\}^{2}}{\alpha_{0,p}(\bm{\psi}_{2} )^{2}+(\lambda_{2}^{2}/4)\left\{ \alpha_{0,p-1}(\bm{\psi}_{2} )-\alpha_{0,p+1}(\bm{\psi}_{2} )\right\}^{2}}. \label{ratio_rhos}
\end{align}
By condition (\ref{cond(iii)}), there exists a function $e_{p}(\bm{\psi} )$ of integer $p$ and positive constant $M>0$ such that $\sup_{p\in \mathbb{N}}|e_{p}(\bm{\psi} )|\leq M$ and 
\begin{align*}
\frac{\alpha_{0,p-1}(\bm{\psi} )-\alpha_{0,p+1}(\bm{\psi} )}{\alpha_{0,p}(\bm{\psi} )}=e_{p}(\bm{\psi} )p^{c},
\end{align*}
where the constant $c$ is given in condition (iii) of Theorem \ref{main_theorem2}. 

Hence, equation (\ref{ratio_rhos}) is expressed as
\begin{align}
\frac{\alpha_{0,p}(\bm{\psi}_{1} )^{2}\left\{ 1+(\lambda_{1}^{2}/4)e_{p}(\bm{\psi}_{1})^{2}p^{2c} \right\}}{\alpha_{0,p}(\bm{\psi}_{2} )^{2}\left\{ 1+(\lambda_{2}^{2}/4)e_{p}(\bm{\psi}_{2})^{2}p^{2c} \right\}}.\label{ratio_rhos2}
\end{align}
Because the domain of $\phi_{1}(p|\bm{\gamma})$ is $S_{1}(\bm{\gamma})=\mathbb{Z}$, and the closure of $S_{1}(\bm{\gamma})$ becomes the extended rational number set $\overline{\mathbb{Z}}$ under the metric $d(x,y)$ described in Section \ref{sec:Main}. Therefore, $p\in\mathbb{Z}$ is allowed to increase to infinity. 

If condition $\{ \alpha_{0,p}(\bm{\psi}_{1})/\alpha_{0,p}(\bm{\psi}_{2})\} p^{c}\rightarrow 0$ $(p\to \infty )$ in (\ref{cond(iii)-2}) holds, then the ratio (\ref{ratio_rhos2}) is bounded above by 
\begin{align}
\frac{\alpha_{0,p}(\bm{\psi}_{1} )^{2}}{\alpha_{0,p}(\bm{\psi}_{2} )^{2}}\left\{ 1+ \frac{M^{2}}{4}p^{2c}\right\} &=\left(\frac{\alpha_{0,p}(\bm{\psi}_{1})}{\alpha_{0,p}(\bm{\psi}_{2} )}p^{c}\right)^{2}\left\{ p^{-2c}+ \frac{M^{2}}{4}\right\} \notag \\
 &\rightarrow 0\quad (p\to\infty).
\end{align}

On the other hand, if $\frac{\alpha_{0,p}(\bm{\psi}_{1})}{\alpha_{0,p}(\bm{\psi}_{2})}p^{-c}\rightarrow \infty$ $(p\to \infty )$, then the ratio (\ref{ratio_rhos2}) is bounded below by 
\begin{align}
\frac{\alpha_{0,p}(\bm{\psi}_{1} )^{2}}{\alpha_{0,p}(\bm{\psi}_{2} )^{2}}\frac{1}{ 1+ (M^{2}/4)p^{2c}} &=\left(\frac{\alpha_{0,p}(\bm{\psi}_{1}}{\alpha_{0,p}(\bm{\psi}_{2} )}p^{-c}\right)^{2} \frac{1}{ p^{-2c}+ (M^{2}/4)}\notag \\
 &\rightarrow \infty \quad (p\to\infty).
\end{align}
Accordingly, under assumption $\bm{\psi}_{1}\ne\bm{\psi}_{2}$, the ratio $\phi_{1}(p|\bm{\gamma}_{1})/\phi_{1}(p|\bm{\gamma}_{2}) \not\rightarrow 1$ $(p\to\infty )$. 

Step 2. Next, we consider two parameter vectors $\bm{\gamma}_{1}$ and $\bm{\gamma}_{2}$ with $\bm{\psi}_{1}=\bm{\psi}_{2}=\bm{\psi}$ and $\lambda_{1}\ne\lambda_{2}$, and set $\phi_{2}(p|\bm{\gamma})=\beta_{p}(\bm{\gamma})$, which is defined in equation (\ref{sine_cosine_moments}). 
Then, by Lemma \ref{lemma1}, for possibly different parameters $\mu_{1}$ and $\mu_{2}$ in $[0,2\pi )$, there exists a sequence $\{ p_{n}\}_{n\in\mathbb{N}}$ with $p_{n}\in\mathbb{N}$ such that $\lim_{n\to\infty}p_{n}=\infty$ and $\lim_{n\to\infty}|p_{n}\mu_{i} \pmod{2\pi}|=0$ $(i=1,2)$. 
By using the sequence $\{ p_{n}\}_{n\in\mathbb{N}}$, it follows that
\begin{align}
\left| \sin (p_{n}\mu_{i}) \frac{\alpha_{0,p_{n}}(\bm{\psi}_{i})}{\alpha_{0,p_{n}-1}(\bm{\psi}_{i} )-\alpha_{0,p_{n}+1}(\bm{\psi}_{i} )}\right| &\leq |\sin (p_{n}\mu_{i})|\frac{1}{\inf_{p\in\mathbb{N}}\left| \frac{\alpha_{0,p-1}(\bm{\psi}_{i} )-\alpha_{0,p+1}(\bm{\psi}_{i} )}{\alpha_{0,p}(\bm{\psi}_{i})}\right|}\notag \\
&\leq M|\sin (p_{n}\mu_{i})|\notag \\
&\rightarrow 0\quad (n\to\infty ). \label{proof_step2}
\end{align}

If $\lambda_{2}\ne 0$, from equation (\ref{proof_step2}), the ratio of $\phi_{2}(p|\bm{\gamma}_{1})$ to $\phi_{2}(p|\bm{\gamma}_{2})$ with $p$ replaced by $p_{n}$ is 
\begin{align*}
\frac{\phi_{2}(p_{n}|\bm{\gamma}_{1})}{\phi_{2}(p_{n}|\bm{\gamma}_{2})}&=\frac{\sin (p_{n}\mu_{1} )\alpha_{0,p_{n}}(\bm{\psi} )+\cos (p_{n}\mu_{1} )\lambda_{1} \{ \alpha_{0,p_{n}-1}(\bm{\psi} )-\alpha_{0,p_{n}+1}(\bm{\psi} )\} /2 }{\sin (p_{n}\mu_{2} )\alpha_{0,p_{n}}(\bm{\psi} )+\cos (p_{n}\mu_{2} )\lambda_{2} \{ \alpha_{0,p_{n}-1}(\bm{\psi} )-\alpha_{0,p_{n}+1}(\bm{\psi} )\} /2 }\notag \\
 &=\frac{\sin (p_{n}\mu_{1} )\cfrac{\alpha_{0,p_{n}}(\bm{\psi} )}{\alpha_{0,p_{n}-1}(\bm{\psi} )-\alpha_{0,p_{n}+1}(\bm{\psi} )}+\cos (p_{n}\mu_{1} )\lambda_{1}/2 }{\sin (p_{n}\mu_{2} )\cfrac{\alpha_{0,p_{n}}(\bm{\psi} )}{\alpha_{0,p_{n}-1}(\bm{\psi} )-\alpha_{0,p_{n}+1}(\bm{\psi} )}+\cos (p_{n}\mu_{2} )\lambda_{2} /2 }\notag \\
 &\rightarrow \frac{\lambda_{1}}{\lambda_{2}}\ne 1\qquad (n\to\infty ).
\end{align*}
If $\lambda_{2}=0$, then $|\phi_{2}(p_{n}|\bm{\gamma}_{1})/\phi_{2}(p_{n}|\bm{\gamma}_{2})|\to \infty$ as $n\to\infty$. 

Step 3. We consider two parameter vectors $\bm{\gamma}_{1}$ and $\bm{\gamma}_{2}$ with $\bm{\psi}_{1}=\bm{\psi}_{2}=\bm{\psi}$, $\lambda_{1}=\lambda_{2}=\lambda$ and $\mu_{1}\ne \mu_{2}$. 
Let $\mathrm{i}$ be the imaginary unit defined as $\mathrm{i}^{2}=-1$, and set $\phi_{3}(\bm{\gamma})=\alpha_{1}(\bm{\gamma})+\mathrm{i}\beta_{1}(\bm{\gamma})=\exp (\mathrm{i}\mu )\{ \alpha_{0,1}(\bm{\psi} )+\mathrm{i}\lambda M_{1}(\bm{\psi} )\}$ which is the polar form of the mean direction where $M_{1}(\bm{\psi} )=\{ \alpha_{0,0}(\bm{\psi} )-\alpha_{0,2}(\bm{\psi} )\} /2$. 
Then, we have 
\begin{align}
&\phi_{3}(\bm{\gamma}_{1})-\phi_{3}(\bm{\gamma}_{2})=\left\{ \exp (\mathrm{i}\mu_{1})-\exp (\mathrm{i}\mu_{2}) \right\}\{ \alpha_{0,1}(\bm{\psi} )+\mathrm{i}\lambda M_{1}(\bm{\psi} )\} \label{step3.1}
\end{align}
Because $|\alpha_{0,1}(\bm{\psi} )+\mathrm{i}\lambda M_{1}(\bm{\psi} )|\geq |\alpha_{0,1}(\bm{\psi} )|>0$ and $\exp (\mathrm{i}\mu_{1})-\exp (\mathrm{i}\mu_{2})\ne 0$, we have $\phi_{3}(\bm{\gamma}_{1})\ne \phi_{3}(\bm{\gamma}_{2})$, which completes the proof.\hfill $\square$

\hspace*{-1em}{\sc Proof of Proposition \ref{thm1}}. 
Condition (i) is obvious. 
It follows from $\alpha_{0,p}(\rho )=E_{0,\rho}\{ \cos (p\Theta )\} =\rho^{p}$ that 
\begin{align}
\frac{\alpha_{0,p-1}(\rho )-\alpha_{0,p+1}(\rho )}{\alpha_{0,p}(\rho )}=\frac{1}{\rho}-\rho .
\end{align}
Thus, conditions (ii) and (iii) hold with $c=0$, which completes the proof \hfill $\square$

To prove Proposition \ref{thm2}, we present a lemma on the modified Bessel function.
\begin{lemma}\label{bessel}
The following results hold:
\begin{align}
&\frac{1}{p!}\left( \frac{\kappa}{2}\right)^{p}\leq I_{p}(\kappa )\leq \frac{1}{p!}\left( \frac{\kappa}{2}\right)^{p}\exp\left( \frac{\kappa^{2}}{4}\right) ,\textrm{ and }\label{ineq01}\\
&\frac{A_{p}(\kappa )}{A_{p}(\kappa^{\prime})}=O\left( \left(\frac{\kappa}{\kappa^{\prime}}\right)^{p}\right) \textrm{ for }\kappa <\kappa^{\prime}\quad  (p\to\infty ), \notag
\end{align}
where $A_{p}(\kappa )=I_{p}(\kappa )/I_{0}(\kappa )$.
\end{lemma}
{\sc Proof}. Using the expansion in \citet[p.349]{MJ09}, we have
\begin{align}
I_{p}(\kappa )&=\sum_{r=0}^{\infty}\dfrac{1}{(p+r)!r!}\left( \frac{\kappa}{2}\right) ^{2r+p} \notag \\
&=\frac{1}{p!}\left( \frac{\kappa}{2}\right) ^{p}\left\{ 1+\frac{1}{(p+1)1!}\left( \frac{\kappa}{2}\right) ^{2}+\frac{1}{(p+2)(p+1)2!}\left( \frac{\kappa}{2}\right) ^{4}+\cdots \right\} \label{expans01} \\
&\leq \frac{1}{p!}\left( \frac{\kappa}{2}\right) ^{p}\left\{ 1+\frac{1}{1!}\left( \frac{\kappa^{2}}{4}\right) +\frac{1}{2!}\left( \frac{\kappa^{2}}{4}\right)^{2}+\frac{1}{3!}\left( \frac{\kappa^{2}}{4}\right)^{3}+\cdots \right\} \notag\\
&=\frac{1}{p!}\left( \frac{\kappa}{2}\right) ^{p}\exp\left( \frac{\kappa^{2}}{4}\right) ,\notag
\end{align}
where the last equation holds from the Maclaurin expansion for $e^{x}$. 
On the one hand, by equation (\ref{expans01}), we have the left side inequality in (\ref{ineq01}). 
Next, we prove the latter. Let $\kappa <\kappa^{\prime}$. Because
\begin{equation*}
\frac{A_{p}(\kappa )}{A_{p}(\kappa^{\prime})}=\frac{I_{0}(\kappa^{\prime})I_{p}(\kappa)}{I_{0}(\kappa )I_{p}(\kappa^{\prime})}, \label{Bessel_ratio}
\end{equation*}
applying the inequality (\ref{ineq01}) to $I_{p}(\kappa )$ and $I_{p}(\kappa^{\prime} )$ yields the result. \qed\vspace{1ex}

{\sc Proof of Proposition} \ref{thm2}. 
We verify conditions (i)--(iii) in Theorem \ref{main_theorem2} to demonstrate the identifiability of the SSvM distribution. 
Because the $p$th cosine-moment is given by $\alpha_{0,p}:=E_{0,\kappa}\{ \cos (p\Theta )\} =I_{p}(\kappa )/I_{0}(\kappa )$, condition (i) holds. In addition, we have
\begin{align}
\frac{\alpha_{0,p-1}(\kappa )-\alpha_{0,p+1}(\kappa )}{\alpha_{0,p}(\kappa )}&=\frac{I_{p-1}(\kappa )-I_{p+1}(\kappa )}{I_{p}(\kappa )}. \label{proof_propSSvM1}
\end{align}
By using the result $I_{\nu -1}(z )- I_{\nu +1}(z ) =(2\nu /z) I_{\nu}(z )$, which is given in \citet[p.376, 9.6.26.]{AS72}, equation (\ref{proof_propSSvM1}) equals $(2p)/\kappa$. 
Therefore, conditions (ii) and equation (\ref{cond(iii)}) in (iii) hold with $c=1$. 

Finally, we verify equation (\ref{cond(iii)-2}) in condition (iii). It follows from Lemma \ref{bessel} that 
\begin{align}
\frac{\alpha_{0,p}(\kappa_{1})}{\alpha_{0,p}(\kappa_{2})}=\frac{I_{p}(\kappa_{1} )/I_{0}(\kappa_{1} )}{I_{p}(\kappa_{2})/I_{0}(\kappa_{2})}=O\left( \left( \frac{\kappa_{1}}{\kappa_{2}}\right)^{p}\right) .
\end{align}
Thus, if $\kappa_{1}>\kappa_{2}$, then $(\alpha_{0,p}(\kappa_{1})/\alpha_{0,p}(\kappa_{2}))p^{-c}\rightarrow \infty$ $(p\to \infty )$. 
By contrast, if $\kappa_{1}<\kappa_{2}$, then $(\alpha_{0,p}(\kappa_{1})/\alpha_{0,p}(\kappa_{2}))p^{c}\rightarrow 0$ $(p\to \infty )$, which implies condition (iii). 
Hence, the result is proved. \qed\vspace{1ex}

{\sc Proof of Proposition} \ref{thm6}. For different parameters $\bm{\gamma}_{\textsc{WS1}}=(\mu_{1},\rho_{\alpha 1},\bar{\rho}_{1},\xi_{1})^{T}$ and $\bm{\gamma}_{\textsc{WS2}}=(\mu_{2},\rho_{\alpha 2},\bar{\rho}_{2},\xi_{2})^{T}$, we consider the following four steps:
Step 1: $\rho_{\alpha 1}\ne\rho_{\alpha 2}$, 
Step 2: $\rho_{\alpha 1}=\rho_{\alpha 2}$ and $\bar{\rho}_{1}\ne \bar{\rho}_{2}$, Step 3: $\rho_{\alpha 1}=\rho_{\alpha 2}$, $\bar{\rho}_{1}=\bar{\rho}_{2}$, and $\xi_{1}\ne \xi_{2}$. Step 4: $\rho_{\alpha 1}=\rho_{\alpha 2}$, $\bar{\rho}_{1}=\bar{\rho}_{2}$, $\xi_{1}=\xi_{2}$, and $\mu_{1}\ne \mu_{2}$. 

First, we consider Step 1 with $\rho_{\alpha 1}\ne\rho_{\alpha 2}$, and set $\phi_{1} (p|\bm{\gamma}_{\textsc{WS}})=\alpha_{p}(\bm{\gamma}_{\textsc{WS}})$ which is the cosine moment. If $\rho_{\alpha 1}>\rho_{\alpha 2}$, then $\phi_{1} (p|\bm{\gamma}_{\textsc{WS1}})/\phi_{1} (p|\bm{\gamma}_{\textsc{WS2}})$ becomes 
\begin{align}
&\frac{p\bar{\rho}_{1}\rho_{\alpha 1}^{p-1}(1-\rho_{\alpha 1}^{2})\cos (p\mu_{1} +\xi_{1} )+\rho_{\alpha 1}^{p}\cos (p\mu_{1})}{p\bar{\rho}_{2}\rho_{\alpha 2}^{p-1}(1-\rho_{\alpha 2}^{2})\cos (p\mu_{2} +\xi_{2} )+\rho_{\alpha 2}^{p}\cos (p\mu_{2})} \notag \\
=&\frac{p\bar{\rho}_{1}(\rho_{\alpha 1}/\rho_{\alpha 2})^{p-1}(1-\rho_{\alpha 1}^{2})\cos (p\mu_{1} +\xi_{1} )+(\rho_{\alpha 1}/\rho_{\alpha 2})^{p-1}\rho_{\alpha 1}\cos (p\mu_{1})/p}{\bar{\rho}_{2}(1-\rho_{\alpha 2}^{2})\cos (p\mu_{2} +\xi_{2} )+\rho_{\alpha 2}\cos (p\mu_{2})/p}\notag \\
&\rightarrow \infty \qquad (p\to \infty ). \notag
\end{align}
The case when $\rho_{\alpha 1}<\rho_{\alpha 2}$ is proved similarly. 

Next, we consider Step 2 with $\rho_{\alpha 1}=\rho_{\alpha 2}=\rho_{\alpha}$ and $\bar{\rho}_{1}\ne \bar{\rho}_{2}$, and set $\phi_{2} (p|\bm{\gamma}_{\textsc{WS}})$ to the squared $p$th mean resultant length, that is $\phi_{2} (p|\bm{\gamma}_{\textsc{WS}})=\rho_{p}(\bm{\gamma}_{\textsc{WS}})^{2}$. 
When $\bar{\rho}_{1}>\bar{\rho}_{2}$, then 
\begin{align}
\frac{\phi_{2} (p|\bm{\gamma}_{\textsc{WS1}})}{\phi_{2} (p|\bm{\gamma}_{\textsc{WS2}})}&=\frac{p^{2}\bar{\rho}_{1}^{2}\rho_{\alpha}^{2(p-1)}(1-\rho_{\alpha}^{2})^{2}+\rho_{\alpha}^{2p}+2p\bar{\rho}_{1}\rho_{\alpha}^{2(p-1)+1}(1-\rho_{\alpha}^{2})\cos \xi_{1}}{p^{2}\bar{\rho}_{2}^{2}\rho_{\alpha}^{2(p-1)}(1-\rho_{\alpha}^{2})^{2}+\rho_{\alpha}^{2p}+2p\bar{\rho}_{2}\rho_{\alpha}^{2(p-1)+1}(1-\rho_{\alpha}^{2})\cos \xi_{2}}\notag \\
&=\frac{\bar{\rho}_{1}^{2}(1-\rho_{\alpha}^{2})^{2}+\rho_{\alpha}^{2}/p^{2}+2\bar{\rho}_{1}\rho_{\alpha}(1-\rho_{\alpha}^{2})(\cos \xi_{1})/p}{\bar{\rho}_{2}^{2}(1-\rho_{\alpha}^{2})^{2}+\rho_{\alpha}^{2}/p^{2}+2\bar{\rho}_{2}\rho_{\alpha}(1-\rho_{\alpha}^{2}) (\cos \xi_{2})/p}\notag \\
&\rightarrow \frac{\bar{\rho}_{1}^{2}}{\bar{\rho}_{2}^{2}}>1 \qquad (p\to\infty ). \notag
\end{align}
The case when $\bar{\rho}_{1}<\bar{\rho}_{2}$ can also be verified similarly. 

Next, we consider Step 3 with $\rho_{\alpha 1}=\rho_{\alpha 2}=\rho_{\alpha}$, $\bar{\rho}_{1}=\bar{\rho}_{2}=\bar{\rho}$, $\xi_{1}\ne\xi_{2}$. $\phi_{3} (p|\bm{\gamma}_{\textsc{WS}})$ is set to the characteristic function $\alpha_{p}(\bm{\gamma}_{\textsc{WS}})+\mathrm{i}\beta_{p}(\bm{\gamma}_{\textsc{WS}})$. For simplicity of notation, let $B_{p}=p\bar{\rho}\rho_{\alpha}^{p-1}(1-\rho_{\alpha}^{2})$. 
If $\xi_{1}\ne \xi_{2}$, then we have
\begin{align}
&\frac{\alpha_{p}(\bm{\gamma}_{\textsc{WS1}})}{\alpha_{p}(\bm{\gamma}_{\textsc{WS2}})}\notag \\
=&\frac{\cos (p\mu_{1} +\xi_{1} )+(\rho_{\alpha}^{p}/B_{p})\cos (p\mu_{1})+\mathrm{i}\left\{ \sin (p\mu_{1} +\xi_{1} )+(\rho_{\alpha}^{p}/B_{p})\sin (p\mu_{1})\right\}}{\cos (p\mu_{2} +\xi_{2} )+(\rho_{\alpha}^{p}/B_{p})\cos (p\mu_{2})+\mathrm{i}\left\{ \sin (p\mu_{2} +\xi_{2} )+(\rho_{\alpha}^{p}/B_{p})\sin (p\mu_{2})\right\}}
 \label{MC_step3}
\end{align}
By Lemma \ref{lemma1}, there exists a sequence $\{ p_{n}\}$ of natural numbers such that $p_{n}\to \infty$ and $|p\mu_{i} \textrm{\upshape (mod 2$\pi$)} |\to 0$ $(i=1,2)$ as $n\to \infty$. Hence, replacing $p$ with $p_{n}$ in equation (\ref{MC_step3}) and taking the limit leads to 
\begin{align}
\frac{\phi_{3} (p_{n}|\bm{\gamma}_{\textsc{WS1}})}{\phi_{3} (p_{n}|\bm{\gamma}_{\textsc{WS2}})} \rightarrow \frac{\cos (\xi_{1})+\mathrm{i}\sin (\xi_{1})}{\cos (\xi_{2})+\mathrm{i}\sin (\xi_{2})} \ne 1 ,\qquad (n\to \infty ), \notag
\end{align}
for which the condition is verified under Step 3. 

Finally, we consider Step 4 with $\rho_{\alpha 1}=\rho_{\alpha 2}=\rho_{\alpha}$, $\bar{\rho}_{1}=\bar{\rho}_{2}=\bar{\rho}$, $\xi_{1}=\xi_{2}=\xi$ and $\mu_{1}\ne \mu_{2}$. Here, we set $\phi_{4} (p|\bm{\gamma}_{\textsc{WS}})=\phi_{3} (p|\bm{\gamma}_{\textsc{WS}})$. Then,
\begin{align*}
\phi_{4} (p|\bm{\gamma}_{\textsc{WS}})&=B_{p}\cos (p\mu +\xi )+\rho_{\alpha}^{p}\cos (p\mu )+\mathrm{i}\left\{ B_{p}\sin (p\mu +\xi )+\rho_{\alpha}^{p}\sin (p\mu ) \right\} \\
 &=B_{p}\exp \{ \mathrm{i}(p\mu +\xi )\} +\rho_{\alpha}^{p}\exp (\mathrm{i}p\mu ).
\end{align*}
Hence, we have
\begin{align}
&\phi_{4} (p|\bm{\gamma}_{\textsc{WS1}})-\phi_{4} (p|\bm{\gamma}_{\textsc{WS2}})\notag \\
&=B_{p}\left\{ \exp (\mathrm{i}(p\mu_{1}+\xi ))-\exp (\mathrm{i}(p\mu_{2}+\xi ))\right\} +\rho_{\alpha}^{p}\left( \exp (\mathrm{i}p\mu_{1})-\exp (\mathrm{i}p\mu_{2})\right) \notag \\
 &=\left( B_{p}\exp (\mathrm{i}\xi)+\rho_{\alpha}^{p}\right) \left\{  \exp (\mathrm{i}p\mu_{1})-\exp (\mathrm{i}p\mu_{2})\right\} .\notag
\end{align}
Because $|B_{p}\exp (\mathrm{i}\xi)+\rho_{\alpha}^{p}|^{2}=B_{p}^{2}+2B_{p}\rho_{\alpha}^{p}\cos \xi +\rho_{\alpha}^{2p}\geq (B_{p}-\rho_{\alpha}^{p})^{2}$ and $B_{p}-\rho_{\alpha}^{p}>0$ for some large $p$, we have $\phi_{4} (p|\bm{\gamma}_{\textsc{WS1}})-\phi_{4} (p|\bm{\gamma}_{\textsc{WS2}})\ne 0$ for some large $p$, which the condition is verified under Step 4. 
Therefore, the family of distributions (\ref{MC01}) is identifiable under $\Gamma_{\textsc{WS}}$. \qed\vspace{1ex}


{\sc Proof of Proposition \ref{thm4}.} Because the marginal distribution of $\Theta$ is the SSWC, from Proposition \ref{thm1} and Remark \ref{remark2}, it is identifiable with respect to the parameter vector $(\mu ,\kappa ,\lambda )^{T}$. 
Hence, we begin with the case of $(\mu_{1},\kappa_{1},\lambda_{1})^{T}=(\mu_{2},\kappa_{2},\lambda_{2})^{T}(=(\mu ,\kappa ,\lambda )^{T})$. 
Set $g(\alpha ,\beta )=\beta (1-\tanh (\kappa)\cos(\theta -\mu))^{1/\alpha}$. 
It follows from the identifiability of Weibull distribution that $f_{\textsc{AL}}(x|\theta ;\alpha_{1},\beta_{1},\kappa )=f_{\textsc{AL}}(x|\theta ;\alpha_{2},\beta_{2},\kappa )$ implies $\alpha_{1}=\alpha_{2}$ and $g(\alpha_{1},\beta_{1})=g(\alpha_{2},\beta_{2})$, which leads to $\beta_{1}=\beta_{2}$. 
This means that for any $(\alpha_{1},\beta_{1})^{T}\ne (\alpha_{2},\beta_{2})^{T}$, $f_{\textsc{AL}}(x|\theta ;\alpha_{1},\beta_{1},\kappa )\ne f_{\textsc{AL}}(x|\theta ;\alpha_{2},\beta_{2},\kappa )$, which completes the proof. \qed\vspace{1ex}

{\sc Proof of Proposition \ref{thm5}.} Because the marginal distribution for $\Theta$ in (\ref{gpar}) is the SSWC distribution, the identifiability with respect to the parameter vector $(\mu, \kappa, \lambda)^T$ holds from Proposition \ref{thm1} and Remark \ref{remark2}. The conditional distribution of the linear part of (\ref{gpar}) is given in equation \eqref{gpar_cond}. We need to check the case with $(\mu_1, \kappa_1 )^T=(\mu_2, \kappa_2)^T=(\mu, \kappa)^T$ and $\delta_1 \ne \delta_2$. First, we assume $\delta_2 > \delta_1$. In this case, the ratio on the conditional distribution of the linear part becomes
\begin{align*}
&\frac{f_{X|\Theta}(x | \theta; \delta_1, \tau_1, \sigma_1)}{f_{X|\Theta}(x | \theta; \delta_2, \tau_2, \sigma_2)}\\
=&
\cfrac{ \frac{1}{\sigma_1 \delta_1} \left( \frac{x}{\sigma_1}\right)^{1/\delta_1 -1} 
\left\{ 
1-\kappa \cos(\theta-\mu)
\right\}
\left[1 + \frac{\tau_1}{\delta_1} \left(\frac{x}{\sigma_1}\right)^{1/\delta_1} \{1 - \kappa \cos(\theta - \mu) \}\right]^{-(\delta_1/\tau_1 + 1)}
}{
 \frac{1}{\sigma_2 \delta_2} \left( \frac{x}{\sigma_2}\right)^{1/\delta_2 -1}
 \left\{ 
1-\kappa \cos(\theta-\mu)
\right\}
\left[1 + \frac{\tau_2}{\delta_2} \left(\frac{x}{\sigma_2}\right)^{1/\delta_2} \{1 - \kappa \cos(\theta - \mu) \}\right]^{-(\delta_2/\tau_2 + 1)}
} \\
=& C_1 x^{1/\delta_1 - 1/\delta_2}
\frac{
\left[1 + \frac{\tau_2}{\delta_2} \left(\frac{x}{\sigma_2}\right)^{1/\delta_2} \{1 - \kappa \cos(\theta - \mu) \}\right]^{(\delta_2/\tau_2 + 1)}
}
{\left[1 + \frac{\tau_1}{\delta_1} \left(\frac{x}{\sigma_1}\right)^{1/\delta_1} \{1 - \kappa \cos(\theta - \mu) \}\right]^{(\delta_1/\tau_1 + 1)}} 
=: K_{1},
\end{align*}
where $C_1 = \sigma_2 ^{1/\delta_2}\delta_2 / (\sigma_1 ^{1/\delta_1}\delta_1 )$. Since $1/\delta_1 - 1/\delta_2 >0$, we observe that $K_{1} \to 0$ as $x \to 0$. We can similarly observe that for the case with $\delta_2 < \delta_1$,  $K_{1} \to \infty$ as $x \to 0$. 
For the next step, we consider the case with $(\mu_1 ,\kappa_1, \delta_1)^T = (\mu_2,\kappa_2, \delta_2)^T = (\mu ,\kappa, \delta)^T$
 and $\tau_1 \ne \tau_2$. 
To simplify expressions, we let $1/\tau_2 -1/\tau_1  =\Delta_{\tau}$ and $C_2(\theta )=1 - \kappa \cos(\theta - \mu)$. 
From the expression $K_{1}$ and the fact that $(\delta /\tau_2)+1 = (\delta /\tau_1) +1+\Delta_{\tau}\delta$, we see that
\begin{align*}
& \frac{f_{X|\Theta}(x | \theta; \delta_1, \tau_1, \sigma_1)}{f_{X|\Theta}(x | \theta; \delta_2, \tau_2, \sigma_2 )}
 =
  C_1
  \frac{
\left[1 + \frac{\tau_2}{\delta} \left(\frac{x}{\sigma_2}\right)^{1/\delta} C_2(\theta)\right]^{(\delta/\tau_2 + 1)}
}
{\left[1 + \frac{\tau_1}{\delta} \left(\frac{x}{\sigma_1}\right)^{1/\delta} C_2(\theta)\right]^{(\delta/\tau_1 + 1)}} \\
=&
C_1 
\left\{
   \frac{
1 + \frac{\tau_2}{\delta} \left(\frac{x}{\sigma_2}\right)^{1/\delta}C_2(\theta) }
{
1 + \frac{\tau_1}{\delta} \left(\frac{x}{\sigma_1}\right)^{1/\delta}C_2(\theta) 
}
\right\}^{\delta/\tau_1+1}
\left[
1+ \frac{\tau_1}{\delta}\left(
\frac{x}{\sigma_2}
\right)^{1/\delta}
C_2(\theta)
\right]^{\Delta_{\tau}} =: C_{1}\times K_{2} \times K_{3} .
\end{align*}
We see that $K_{2} =O(1)$ as $x\to \infty$. For $\tau_2 < \tau_1$ which indicates $\Delta_{\tau} >0$, we have 
 $K_{3} \to \infty$ as $x \to \infty$. On the contrary, when $\tau_2 > \tau_1$ which indicates $\Delta_{\tau} <0$, we have 
 $K_{3} \to 0$ as $x \to \infty$.
For the final step, we consider the case with $(\mu_1 ,\kappa_1, \delta_1,\tau_1)^T = (\mu_2,\kappa_2, \delta_2,\tau_2)^T = (\mu ,\kappa, \delta ,\tau)^T$ and $\sigma_1 \ne \sigma_2$. We observe that
\begin{align*}
 \frac{f_{X|\Theta}(x | \theta; \delta, \tau, \sigma_1)}{f_{X|\Theta}(x | \theta; \delta, \tau, \sigma_2 ) }
 =
\frac{\sigma_2^{1/\delta}\delta }{\sigma_1^{1/\delta}\delta }
  \frac{
\left[1 + \frac{\tau}{\delta} \left(\frac{x}{\sigma_2}\right)^{1/\delta} C_2(\theta)\right]^{(\delta/\tau + 1)}
}
{\left[1 + \frac{\tau}{\delta} \left(\frac{x}{\sigma_1}\right)^{1/\delta} C_2(\theta)\right]^{(\delta/\tau + 1)}} .
\end{align*}
 Evaluating the ratio of the conditional distributions at $x=0$, we have
 \begin{align*}
  \frac{f_{X|\Theta}(0 | \theta; \delta, \tau, \sigma_1)}{f_{X|\Theta}(0| \theta; \delta, \tau, \sigma_2 ) }
 =
\frac{\sigma_2^{1/\delta}}{\sigma_1^{1/\delta}} \ne 1,
 \end{align*} 
 which completes the proof of Proposition \ref{thm5}. \hfill $\square$

\section{Derivation of $p$th cosine and sine moments}\label{appendix:A}
\hspace*{1em}In this section, we describe how to represent $p$th cosine and sine moments in a circular model with location $\mu$ by using $p$th cosine and sine moments in the circular model with $\mu =0$. 
We consider the density \eqref{sscd} with $\mu =0$,
\begin{align}
f^{*}(\theta |\bm{\gamma}^{*})&:=f_{0}(\theta |\bm{\psi} )\left\{ 1+\lambda \sin (\theta )\right\} ,
\end{align}
where $\bm{\gamma}^{*}=(\bm{\psi}^{T},\lambda )^{T}$. 
In addition, we denote its $p$th cosine and sine moments as $\alpha_{p}^{*}(\bm{\gamma}^{*}):=E^{*}\{\cos (p\Theta)\}$ and $\beta_{p}^{*}(\bm{\gamma}^{*}):=E^{*}\{\sin (p\Theta)\}$, $(p\in \mathbb{Z})$. 
Because the density \eqref{sscd} is written as $f(\theta |\bm{\gamma})=f^{*}(\theta -\mu |\bm{\gamma}^{*})$, the $p$th cosine moment in the density \eqref{sscd} is given by
\begin{align}
\alpha_{p}(\bm{\gamma})&=\int_{\mu-\pi}^{\mu+\pi}\cos (p\theta )f^{*}(\theta -\mu |\bm{\gamma}^{*}) d\theta .\label{ap:eq2}
\end{align}
Letting $\theta -\mu =\eta$, equation \eqref{ap:eq2} becomes
\begin{align}
 &=\int_{-\pi}^{\pi}\cos (p(\mu +\eta ))f^{*}(\eta |\bm{\gamma}^{*})d\eta \notag\\
 &=\int_{-\pi}^{\pi}\left\{ \cos (p\mu )\cos(p\eta )-\sin (p\mu )\sin (p\eta )\right\}f^{*}(\eta |\bm{\gamma}^{*}) d\eta \notag\\
&=\cos (p\mu )E^{*}\left\{ \cos (p\eta )\right\} -\sin (p\mu)E^{*}\left\{ \sin (p\eta )\right\} \notag \\
&=\cos (p\mu )\alpha_{p}^{*}(\bm{\gamma}^{*})-\sin (p\mu )\beta_{p}^{*}(\bm{\gamma}^{*}).
\end{align}

Similarly, the $p$th sine moment is 
\begin{align}
\beta_{p}(\bm{\gamma})&=\int_{-\pi}^{\pi}\sin (p\mu +p\eta )f^{*}(\eta |\bm{\gamma}^{*})d\eta \notag \\
              &=\int_{-\pi}^{\pi}\left\{ \sin (p\mu )\cos (p\eta )+\cos (p\mu )\sin (p\eta )\right\} f^{*}(\eta |\bm{\gamma}^{*})d\eta \notag \\
              &=\sin (p\mu )E^{*}\left\{ \cos (p\eta )\right\} +\cos (p\mu)E^{*}\left\{ \sin (p\eta )\right\} \notag \\
              &=\sin (p\mu )\alpha_{p}^{*}(\bm{\gamma}^{*}) +\cos (p\mu)\beta_{p}^{*}(\bm{\gamma}^{*}). \notag 
\end{align}
This result shows that if we know the $p$th cosine and sine moments when $\mu =0$, we can automatically derive the $p$th cosine and sine moments of the circular model with location $\mu$.

\end{document}